\documentclass[11pt,
]{amsart}
\usepackage{
amssymb, graphicx
}
\usepackage{array}
\usepackage{graphicx,color}
\usepackage{bm}
\newtheorem{thm}{Theorem}

\newtheorem{remark}{Remark}
 
\date{\today}

\newcommand{\be}[1]{\begin{equation}\label{#1}} 
\newcommand{\ee}{\end{equation}} 

\title[Uniqueness of a transversely isotropic perturbation]
{Unique determination of a transversely isotropic perturbation in a linearized inverse boundary value problem for elasticity}
\author[Y. Yang]{Yang Yang}
\address{Department of Computational Mathematics, Science and Engineering, Michigan State University, East Lansing,
    MI 48824, USA (\tt{yangy5@msu.edu}).}
  \author[J. Zhai]{Jian Zhai}
\address{Department of Mathematics,
 University of Washington, Seattle, WA 98155, USA
  (\tt{jian.zhai@outlook.com}).}

\thanks{The research of YY was partly supported by NSF Grant DMS-1715178, AMS-Simons travel grant, and start-up fund from Michigan State University.}

\begin{document}

\maketitle
\begin{abstract}
We consider a linearized inverse boundary value problem for the elasticity system. From the linearized Dirichlet-to-Neumann map at zero frequency, we show that a transversely isotropic perturbation of a homogeneous isotropic elastic tensor can be uniquely determined. From the linearized Dirichlet-to-Neumann map at two distinct positive frequencies, we show that a transversely isotropic perturbation of a homogeneous isotropic density can be identified at the same time.

\end{abstract}

\section{Introduction and Main result}


In this paper, we investigate the problem of determining interior material property of an elastic body from boundary measurements.
We will consider the linearized inverse boundary value problem for the equation
$$
\text{div}(\mathbf{C}\nabla u)+\omega^2\rho u =0,
$$
which reads in components as 
\begin{equation} \label{eq:elastic1}
\partial_j C_{ijkl}(x) \partial_k u_l(x) +\omega^2\rho_{ik}(x)u_k(x)= 0, \quad\quad i = 1,2,3.
\end{equation}
Here $\omega \geq 0$ is the frequency, $u$ is the displacement vector, $\rho = (\rho_{ik})$ is a symmetric matrix representing the density of mass; $\mathbf{C}=(C_{ijkl}))$ is the elastic tensor whose components obey the symmetry conditions
\begin{equation} \label{eq:Csym}
C_{ijkl} = C_{jikl} = C_{klij}.
\end{equation}
We have used Einstein's summation convention in \eqref{eq:elastic1} such that repeated indices are summed up over $\{1,2,3\}$. Note that $\mathbf{C}$ with the above symmetry has a total number of 21 linearly independent components. For a fixed $\omega$, the case $\omega=0$ corresponds to the governing equations for linear elasticity in equilibrium, while the case $\omega>0$ represents the time-harmonic elastic wave with frequency $\omega$.

Let $\Omega$ be an open bounded domain in $\mathbb{R}^3$ with $C^{1,1}$ boundary $\partial \Omega$. Suppose the density of mass $\rho$ and the elastic tensor $\mathbf{C}$ are both bounded, in the sense that $\rho_{ik}, \, C_{ijkl}\in L^\infty(\Omega)$ for all $i,j,k,l \in \{1,2,3\}$.
We further assume that the density of mass $\rho$ and the elasticity tensor $\mathbf{C}$ satisfy the following positivity conditions: 
there exists $\delta>0$ such that for any real-valued $3$-vector $\sigma = (\sigma_1,\sigma_2,\sigma_3)$,
$$
   \sum_{i,k=1}^3 \rho_{ik} \sigma_{i}\sigma_{k} \geq \delta \sum_{i=1}^3 \sigma_{i}^2;
$$
and for any 
$3\times 3$ real-valued symmetric matrix $(\varepsilon_{ij})$, 
\[
   \sum_{i,j,k,l=1}^3C_{ijkl}\varepsilon_{ij}\varepsilon_{kl}\geq\delta\sum_{i,j=1}^3\varepsilon_{ij}^2.
\]
If $\omega^2$ is not a Dirichlet eigenvalue of the operator $-\rho^{-1} \text{div}(\mathbf{C}\nabla \cdot)$, then for any $f\in H^{1/2}(\partial\Omega)$, standard elliptic theory ensures a unique solution $u^f\in H^{1}(\Omega)$ to the boundary value problem
$$ \left\{
\begin{array}{rcl}
\partial_j C_{ijkl}(x) \partial_k u^f_l(x) + \omega^2 \rho_{ik} (x) u^f_k(x) & = & 0\text{~~~in~}\Omega, \quad\quad i = 1,2,3 \vspace{1ex}\\
u^f |_{\partial\Omega} & = & f.
\end{array}
\right. $$
We define the Dirichlet-to-Neumann map (DN map) $\Lambda_{\mathbf{C}}$ by
$$
\Lambda_{\mathbf{C},\rho,\omega}: f\mapsto C_{ijkl}\nu_j\partial_k u^f_l\vert_{\partial\Omega}
$$
where $\nu=(\nu_1,\nu_2,\nu_3)$ denotes the outer unit normal vector to $\partial\Omega$. It follows that $\Lambda_{\mathbf{C},\rho,\omega}: H^{1/2}(\partial\Omega)\rightarrow H^{-1/2}(\partial\Omega)$ is a bounded linear operator, and the equivalent weak formulation is
$$
\langle \Lambda_{\mathbf{C},\rho,\omega} f,g \rangle = \int_{\Omega} C_{ijkl}\partial_i u^f_j \partial_k u^g_l - \omega^2\rho_{ik} u^f_i u^g_k\,dx.
$$
for any $f,g \in H^{1/2}(\partial\Omega)$. We are interested in determining $\mathbf{C},\rho$ from $\Lambda_{\mathbf{C},\rho,\omega}$
This is related to the invertibility of the non-linear map $(\mathbf{C},\rho) \mapsto \Lambda_{\mathbf{C},\rho,\omega}$. The question is difficult in the general setting, so it is commonly studied under additional a-priori information.

\bigskip
\textbf{The case $\boldsymbol{\omega=0}$.}
Note that when $\omega=0$, the density $\rho$ does not appear in the equation \eqref{eq:elastic1}, thus one can only expect to recover information on $\mathbf{C}$. We henceforth write $\Lambda_{\mathbf{C},\rho,0}$ as $\Lambda_{\mathbf{C}}$ for the ease of notation.

 We say the elastic tensor $\mathbf{C}$ (or the medium) is \textit{homogeneous} if it is a constant tensor (that is, independent of $x$); it is \textit{isotropic} if it can be written as
$$
C_{ijkl}(x) = \lambda(x) \delta_{ij}\delta_{kl} + \mu(x)(\delta_{ik}\delta_{jl} + \delta_{il}\delta_{jk})
$$
where the two functions $\lambda(x)$ and $\mu(x)$ are known as \textit{Lam\'{e} parameters}; and it is \textit{fully anisotropic} if the components $C_{ijkl}$ are subject to no other relations other than \eqref{eq:Csym}. For isotropic $\mathbf{C}$, a global uniqueness result can be found \cite{IY} in dimension two. The problem remains open in dimension three, yet some special cases have been tackled. Among them, 
Nakamura and Uhlmann \cite{NU2} proved uniqueness when the Lam\'{e}
parameters are smooth and $\mu(x)$ is close to a positive
constant, see \cite{ER} for a similar result by Eskin and Ralston and \cite{IUY} for a partial data version; uniqueness for recovering piecewise constant Lam\'{e} parameters was proved in \cite{BFMRV,BFV}; and some boundary determination results were shown in \cite{LN,NU3, NTU} .
For fully anisotropic $\mathbf{C}$, uniqueness was proved in \cite{CHN} for piecewise homogeneous medium. 

It is widely believed that a fully anisotropic $\mathbf{C}$ without additional assumption cannot be uniquely recovered. 
For the inverse conductivity problem, that is the problem to determine the coefficients
$\gamma=(\gamma_{ij}(x))$ in the equation
$$
\partial_i(\gamma_{ij}(x)\partial_ju(x))=0
$$
from the associated Dirichlet-to-Neumann map, it is known that an anisotropic $\gamma(x)$ can at best be determined up to boundary-fixing diffeomorphisms \cite{greenleaf}. 
In contrast, many anisotropic elastic materials have extra structural symmetries which cannot be preserved under diffeomorphisms.
It is therefore important to study the uniqueness of elasticity parameters with extra symmetries in anisotropy. We list some frequently considered anisotropies with symmetries in the table below, see \cite[Chapter 2.6]{Bauchau2009book} \cite[Chapter 3.4]{Reddy2002book} for detailed description. It is worth mentioning that these concepts of anisotropy are purely Cartesian (in a prescribed coordinate system $(x_1, x_2, x_3)$).
\begin{center}
  \begin{tabular}{ | c | m{8cm} | m{3.5cm} |}
    \hline
    Type of anisotropy & \hspace{3cm} Symmetry & Number of independent components \\ \hline
    isotropic & radial symmetry & \hspace{1.5cm} 2 \\ \hline
    cubic & three mutually orthogonal planes of reflection
symmetry plus $\frac{\pi}{2}$ rotation symmetry with respect
to those planes & \hspace{1.5cm} 3 \\ \hline
    transversely isotropic & three mutually orthogonal planes of reflection
symmetry and one symmetry axis perpendicular to one symmetry plane& \hspace{1.5cm} 5 \\ \hline
  orthotropic (orthorhombic)   & three mutually orthogonal planes of reflection
symmetry & \hspace{1.5cm} 9 \\ \hline
  monoclinic   & one plane of reflection symmetry & \hspace{1.4cm} 13 \\ \hline
  fully anisotropic   & no symmetry & \hspace{1.4cm} 21 \\ 
  \hline
  \end{tabular}
\end{center}

In this article, we investigate the linearization of the map $\mathbf{C} \mapsto \Lambda_{\mathbf{C}}$ at a homogeneous isotropic elastic tensor. More specifically, suppose
$$
\mathbf{C}(x)=\mathbf{C}^0+\delta\mathbf{C}(x)
$$
where $\mathbf{C}^0= \lambda^0 \delta_{ij}\delta_{kl} + \mu^0(\delta_{ik}\delta_{jl} + \delta_{il}\delta_{jk})$ is a homogeneous, isotropic background tensor with Lam\'{e} parameters $(\lambda^0, \mu^0)$ satisfying
\begin{equation}\label{cond}
\mu^0>0,\quad\quad 3\lambda^0+2\mu^0>0,
\end{equation}
 and $\delta\mathbf{C}(x)$ is viewed as a perturbation term with components $\delta C_{ijkl}(x)$. It is routine to verify that the map $\mathbf{C} \mapsto \Lambda_{\mathbf{C}}$ is Frech\'{e}t differentiable at $\mathbf{C}^0$ (we refer to \cite{Ike1} for more details), and the Frech\'{e}t derivative
$$
\dot{\Lambda}_{\mathbf{C}^0}: L^\infty(\Omega)\ni\delta\mathbf{C} \mapsto \dot{\Lambda}_{\mathbf{C}^0}(\delta\mathbf{C})\in \mathcal{L}(H^{1/2}(\partial\Omega),H^{-1/2}(\partial\Omega))
$$
is characterized by
\begin{equation} \label{eq:intiden1}
\langle\dot{\Lambda}_{\mathbf{C}^0}(\delta\mathbf{C})f,g\rangle=\int_\Omega \delta C_{ijkl}(x) \, \partial_i u_j(x) \, \partial_k v_l(x) \,dx
\end{equation}
where $u$ (resp. $v$) solves
\begin{equation} \label{eq:elastic2}
\left\{
\begin{array}{rcl}
\mu^0 \Delta u + (\lambda^0 + \mu^0) \nabla \nabla\cdot u & = & 0 \quad \text{ in } \Omega \\
u|_{\partial\Omega} & = & f.
\end{array}
\right.
\left( \text{ resp. }
\left\{
\begin{array}{rcl}
\mu^0 \Delta v + (\lambda^0 + \mu^0) \nabla \nabla\cdot v & = & 0 \quad \text{ in } \Omega \\
v|_{\partial\Omega} & = & g.
\end{array}
\right.
\right)
\end{equation}
The question we are interested in is whether the linearized map $\dot{\Lambda}_{\mathbf{C}^0}$ is injective on anisotropy perturbations with certain symmetry. 
It was proved in \cite{Ike1} that the linearization $\dot{\Lambda}_{\mathbf{C}^0}$ is injective on isotropic perturbations.  Our main theorem (see Theorem \ref{thm:main} below) generalizes this injectivity result from isotropic perturbtions to transversely isotropic perturbations.




A \textit{transversely isotropic} material is one with physical properties that are symmetric about an axis that is normal to a plane of isotropy. It is also known as ``polar anisotropic'' since the material properties are the same in all directions within the transverse plane. Examples of transversely isotropic materials include some piezoelectric materials and fiber-reinforced composites where all fibers are in parallel. Geological layers of rocks are often interpreted as being transversely isotropic as well in terms of their effective properties. Transversely isotropic materials have been extensively studied in geophysical literature, see \cite{TI2,TI3,TI1,TI4} and the references therein.

As is indicated in the above table, transverse isotropy means the elasticity have three mutually orthogonal planes of reflection symmetry and one symmetry axis perpendicular to one of the three symmetry planes. 
Assume the symmetry axis is $x_3$, then
$\delta\mathbf{C}$ obeys the invariance
$$
Q_{ip}Q_{jq}Q_{kr}Q_{ls}\delta C_{pqrs}=\delta C_{ijkl},
$$
where $Q$ can take any of the following reflection and rotation matrices.
\begin{equation}\label{Qmatrix}
\begin{split}
&\left(\begin{array}{ccc}
-1 & 0 &0\\
0  & 1 &0\\
0  & 0 &1
\end{array}\right), \quad
\left(\begin{array}{ccc}
1  & 0 &0\\
0  & -1&0\\
0  & 0 &1
\end{array}\right), \quad
\left(\begin{array}{ccc}
1  & 0 &0\\
0  & 1 &0\\
0  & 0 &-1
\end{array}\right), \\
&\quad\quad\quad\quad\quad\left(\begin{array}{ccc}
\cos\theta & -\sin\theta &0\\
\sin\theta &\cos\theta&0\\
0&0&1
\end{array}\right), \,0\leq \theta \leq 2\pi.
\end{split}
\end{equation}
Writing the above invariance componentwisely results in $9$ non-zero components in $\delta\mathbf{C}$
$$
\delta C_{1111}, \,\delta C_{2222}, \,\delta C_{3333}, \,\delta C_{1122}, \,\delta C_{1133}, \,\delta C_{2233}, \,\delta C_{1212},\,\delta C_{1313}, \,\delta C_{2323}
$$
subject to $4$ linear relations
\begin{equation}
\label{eq:linrel}
\begin{split}
\delta C_{1111} = \delta C_{2222}, \quad\quad &  \delta C_{1133} = \delta C_{2233}, \\
\delta C_{1313}=\delta C_{2323}, \quad\quad & \delta C_{1212} = \frac{1}{2}(\delta C_{1111}-\delta C_{1122}).
\end{split}
\end{equation}
Hence a traversally isotropic $\delta\mathbf{C}$ has only $5$ linearly independent components. We will prove these independent components are uniquely determined by the linearized map $\dot{\Lambda}_{\mathbf{C}^0}$. More precisely, we show
%

\begin{thm} \label{thm:main}
Let $\mathbf{C}^0 = \lambda^0 \delta_{ij}\delta_{kl} + \mu^0(\delta_{ik}\delta_{jl} + \delta_{il}\delta_{jk})$ be homogeneous and isotropic with Lam\'{e} parameters $(\lambda^0, \mu^0)$ satisfying $(\ref{cond})$.
If $\dot{\Lambda}_{\mathbf{C}^0}(\delta\mathbf{C})=0$ and $\delta\mathbf{C} \in L^\infty(\Omega)$ is transversely isotropic with known axis of symmetry, then $\delta\mathbf{C}=0$.
\end{thm}
The $5$ linearly independent components of $\delta\mathbf{C}$ we will determine are $\delta C_{1111}$, $\delta C_{1122}$, $\delta C_{1133}$, $\delta C_{1313}$, and $\delta C_{3333}$.

Injectivity of the linearized map $\dot{\Lambda}_{\mathbf{C}^0}$ has been studied in previous literature. In dimension two or higher, it is known that $\dot{\Lambda}_{\mathbf{C}^0}$ is injective on isotropic $\delta\mathbf{C}$ \cite{Ike1}. Theorem \ref{thm:main} can be viewed as generalization of such injectivity result from isotropic perturbations to transversely isotropic perturbations in dimension three. Note that Theorem \ref{thm:main} has greatly increased the number of independent parameters that can be simultaneously identified by $\dot{\Lambda}_{\mathbf{C}_0}$ -- from $2$ Lam\'{e} parameters in the isotropic case  \cite{Ike1} to $5$ independent parameters in the transversely isotropic case. In dimension two, Ikehata \cite{Ike2,Ike3, Ike4} characterized the injectivity with general anisotropic $\mathbf{C}^0$.

\bigskip
\textbf{The case $\boldsymbol{\omega>0}$.} The time-harmonic case has important application in (reflection) seismology, where one hopes to recover the material parameters of the Earth's subsurface areas from vibroseis data. Unique determination of piecewise homogeneous isotropic parameters from $\Lambda_{\mathbf{C,\rho,\omega}}$ was established in \cite{Beretta2017uniqueness}; unique determination of an anisotropic density with homogeneous isotropic elastic tensor was proved in \cite{Ruiz2018uniqueness}. On the other hand, the inverse boundary value problem for the dynamic elasticity system has been considered in \cite{rachele2000inverse,rachele2003uniqueness,stefanov2017local,bhattacharyya2018local}, 
.

We still consider the linearization of the map $(\mathbf{C},\rho)\mapsto \Lambda_{\mathbf{C,\rho,\omega}}$ at a homogeneous and isotropic $(\mathbf{C}^0, \rho^0)$. Assume $\omega^2$ is not a eigenvalue of $-(\rho^{0})^{-1} \text{div}(\mathbf{C}^0\nabla \cdot)$. The Frech\'{e}t derivative
$$
\dot{\Lambda}_{\mathbf{C}^0,\rho^0,\omega}: L^\infty(\Omega)\ni(\delta\mathbf{C},\delta\rho) \mapsto \dot{\Lambda}_{\mathbf{C}^0,\rho^0,\omega}(\delta\mathbf{C},\delta\rho)\in \mathcal{L}(H^{1/2}(\partial\Omega),H^{-1/2}(\partial\Omega))
$$
is characterized by
\begin{equation} \label{eq:intidenk1}
\langle\dot{\Lambda}_{\mathbf{C}^0,\rho^0,\omega}(\delta\mathbf{C},\delta\rho)f,g\rangle=\int_\Omega \delta C_{ijkl}(x) \, \partial_i u_j(x) \, \partial_k v_l(x) -\omega^2\delta\rho_{ik}u_iv_k\,dx,
\end{equation}
where $u$ and $v$ solve respectively
\begin{equation} \label{eq:elastic2k}
\begin{split}
\left\{
\begin{array}{rcl}
\mu^0 \Delta u + (\lambda^0 + \mu^0) \nabla \nabla\cdot u +\omega^2\rho^0 u& = & 0 \quad \text{ in } \Omega \\
u|_{\partial\Omega} & = & f,
\end{array}
\right.\\
 \text{ }\\
\left\{
\begin{array}{rcl}
\mu^0 \Delta v + (\lambda^0 + \mu^0) \nabla \nabla\cdot v +\omega^2\rho^0 v& = & 0 \quad \text{ in } \Omega \\
v|_{\partial\Omega} & = & g.
\end{array}
\right.
\end{split}
\end{equation}
By definition, a transversely isotropic $\delta\rho$ with symmetry axis $x_3$ has the property
$$
Q_{ip}Q_{jq}\delta\rho_{pq}=\delta\rho_{ij}
$$
for any $Q$ of the forms \eqref{Qmatrix}. Then it
can be written as
\[
\delta\rho=\left(\begin{array}{ccc}\delta\rho_{11}&&\\
&\delta\rho_{11}&\\
&&\delta\rho_{33}
\end{array}\right).
\]

 For the inverse boundary value problem for time-harmonic acoustic wave equation 
\[
\nabla\cdot\gamma\nabla u+\omega^2q u=0,
\]
it is known that the simultaneous recovery of $\gamma$ and $q$ requires two freqeuncy data $\Lambda_{\gamma,q,\omega_1},\,\Lambda_{\gamma,q,\omega_2}$ \cite{nachman1988reconstructions}. We will also use two frequency data in the following. Although single frequency data is enough for recovering piecewise homogeneous parameters \cite{Beretta2017uniqueness}, we do believe a second frequency is necessary for our problem. More precisely, we will prove
\begin{thm} \label{thm:maink}
Let $\mathbf{C}^0 = \lambda^0 \delta_{ij}\delta_{kl} + \mu^0(\delta_{ik}\delta_{jl} + \delta_{il}\delta_{jk})$ be homogeneous and isotropic with Lam\'{e} parameters $(\lambda^0, \mu^0)$ satisfying $(\ref{cond})$, and $\rho^0_{ik}=\rho^0\delta_{ik}$ be a homogeneous isotropic density.
If $\dot{\Lambda}_{\mathbf{C}^0,\rho^0,\omega_i}(\delta\mathbf{C},\delta\rho)=0$ for two distinct positive frequencies $\omega_1,\omega_2$ and $(\delta\mathbf{C},\delta\rho)\in L^\infty(\Omega)$ are transversely isotropic with known axis of symmetry, then $(\delta\mathbf{C},\delta\rho)=0$.
\end{thm}

%
%
%

The proofs of Theorem \ref{thm:main} and Theorem \ref{thm:maink} are based on construction of the \textit{Complex Geometric Optics} (CGO) solutions for the system \eqref{eq:elastic2} and \eqref{eq:elastic2k}, respectively. CGO solutions were initiated by Sylvester and Uhlmann \cite{SU1} in their solving Calder\'{o}n's inverse conductivity problem \cite{Calderon}. Solutions of this type with $\omega=0$ were introduced in \cite{Ike1} for the elasticity system with constant coefficients, and in \cite{ER,NU2} for variable coefficients. Solutions with $\omega> 0$ were utilized in \cite{Ruiz2018uniqueness}.\\

There are notable differences in the construction of CGO solutions for $\omega=0$ and $\omega> 0$. To see this, consider the solution $\psi$ to the scalar wave equation
\[
(\lambda^0+2\mu^0)\Delta\psi+\rho^0 \omega^2\psi=0.
\]
For $\omega=0$, a CGO solution can be taken as $\psi=e^{\zeta\cdot x}$ with $\zeta\in \mathbf{C}^3$, $\zeta\cdot\zeta =0$, then $u:=\nabla\psi$ is a divergence free solution of the elasticity equation in \eqref{eq:elastic2}, which corresponds to an \textit{S}-wave. For $\omega>0$, a similar CGO solution can be constructed as $\psi=e^{\zeta'\cdot x}$ but with $\zeta'\in \mathbf{C}^3$, $\zeta'\cdot\zeta'=\frac{\omega^2\rho^0}{\lambda^0+2\mu^0}$. Then $u:=\nabla\psi$ remains a solution to the elasticity equation in \eqref{eq:elastic2k}, but is not divergence free any more. In fact, this solution $u$ corresponds to a \textit{P}-wave. Of course, this argument is just heuristic as the equation with $\omega=0$ does not really describe waves; but it demonstrates the difference between the cases $\omega=0$ and $\omega>0$. The proofs of Theorem \eqref{thm:main} and Theorem \eqref{thm:maink} are therefore presented separately due to some essential differences in the construction of CGO solutions. The rest of the paper is devoted to these proofs. 

\section{Zero frequency case: Proof of Theorem \ref{thm:main}}

We prove Theorem \ref{thm:main} in this section.
In view of \eqref{eq:intiden1}, $\dot{\Lambda}_{\mathbf{C}^0}(\delta \mathbf{C})=0$ implies
\begin{equation} \label{eq:intiden}
\int_\Omega \delta C_{ijkl} \partial_i u_j \partial_k v_l \,dx=0,
\end{equation}
for any $u,v$ satisfying \eqref{eq:elastic2}. The key ingredient of our proof is constructing CGO solutions to \eqref{eq:elastic2} and inserting them into \eqref{eq:intiden} to obtain sufficiently many linearly independent equations in the $5$ independent components of $\delta\mathbf{C}$. For the ease of notation, we abbreviate $\mathbf{C}$ for $\delta \mathbf{C}(x)$ and $C_{ijkl}$ for $\delta C_{ijkl}(x)$ from now on. We reserve the letter $i$ for an index and write the bold face $\mathbf{i}$ for the imaginary unit.

\bigskip

\noindent \textbf{Step 1.}
Set 
\begin{equation}\label{phase1}
\zeta^{(1)} := \mathbf{i}(s,0,t) + (-t,0,s), \quad\quad\quad 
\zeta^{(2)} := \mathbf{i}(s,0,t) - (-t,0,s),
\end{equation}
and
$$
a^{(1)} = a^{(2)} = a = (0,1,0).
$$ 
We take
$$
u = a^{(1)}e^{\zeta^{(1)} \cdot x},\quad\quad v = a^{(2)}e^{\zeta^{(2)} \cdot x}.
$$
The choice of $\zeta^{(i)}\in\mathbb{C}^3$ ensures $\zeta^{(1)} \cdot \zeta^{(1)} = \zeta^{(2)} \cdot \zeta^{(2)} = 0$, hence $\Delta u = \Delta v = 0$. The choice of $a$ ensures $a\perp \Re \zeta^{(i)}$, $a\perp \Im \zeta^{(i)}$ for $i=1,2$, hence $\nabla\cdot u=\nabla\cdot v=0$. This verifies that $u,v$ defined in this manner satisfy the equations $(\ref{eq:elastic2})$.

Substituting $u$ and $v$ into $(\ref{eq:intiden})$, we have
\begin{align*}
0 = & \int_\Omega C_{ijkl} a_i \zeta^{(1)}_j\; a_k \zeta^{(2)}_l e^{(\zeta^{(1)}+\zeta^{(2)}) \cdot x} \, dx \\
= & \int_\Omega [C_{1212}(\mathbf{i}s-t)(\mathbf{i}s+t)+C_{2323}(\mathbf{i}t+s)(\mathbf{i}t-s)] e^{2\mathbf{i}(s,0,t)\cdot x}\, dx\\
= & \int_\Omega (s^2+t^2)(-C_{1212}-C_{1313})e^{2\mathbf{i}(s,0,t)\cdot x}\, dx.
\end{align*}
This implies the Fourier transform $\mathcal{F}[\chi_\Omega (C_{1212}+C_{1313})](-2(s,0,t))=0$ for any $s,t$ with $s^2+t^2 \neq 0$. Here $\chi_\Omega$ is the characteristic function of the domain $\Omega$. Since $s,t$ can be any real number, this Fourier transform vanishes on the punctured $x_1 x_3$-plane. The axial symmetry with respect to $x_3$-axis in the definition of transversal isotropy allows one to obtain similar vanishing result in any plane containing $x_3$-axis. We conclude $\mathcal{F}[\chi_\Omega (C_{1212}+C_{1313})](\xi)=0$ for any $\xi \neq 0$. This forces
%
%
%
%
$$
C_{1212}+C_{1313}=0 \quad\quad\quad \text{ in } \Omega.
$$
Using the relation $C_{1212} = \frac{1}{2}(C_{1111}-C_{1122})$ in \eqref{eq:linrel} we have
\begin{equation}\label{eq:r5}
C_{1111} - C_{1122} + 2C_{1313} = 0 \quad\quad\quad \text{ in } \Omega.
\end{equation}

\bigskip
\noindent \textbf{Step 2.}
Take
$$
u = \zeta^{(1)} e^{\zeta^{(1)} \cdot x}, \quad\quad\quad 
v =  \zeta^{(2)} e^{\zeta^{(2)} \cdot x},
$$
with $\zeta^{(1)},\,\zeta^{(2)}$ defined in $(\ref{phase1})$. One still has $\Delta u = \Delta v = 0$ as before. On the other hand, the $i$-th component of $u$ (resp. $v$) is $u_i = \zeta^{(1)}_i e^{\zeta^{(1)}\cdot x}$ (resp. $v_i = \zeta^{(2)}_i e^{\zeta^{(2)}\cdot x}$), $i=1,2,3$.
The derivatives of these components are
$$
\partial_i u_j = \zeta^{(1)}_i \zeta^{(1)}_j e^{\zeta^{(1)} \cdot 
x}, \quad\quad\quad 
\partial_k v_l = \zeta^{(2)}_k \zeta^{(2)}_l e^{\zeta^{(2)} \cdot x}.
$$
Then $\nabla\cdot u = \zeta^{(1)} \cdot \zeta^{(1)} e^{\zeta^{(1)} \cdot x} = 0$ and likewise $\nabla\cdot v = 0$. We see that $u,v$ solve \eqref{eq:elastic2}.

Inserting $u,v$ into the integral identity \eqref{eq:intiden}, we obtain
\begin{align*}
0 = \int_\Omega & C_{ijkl} \zeta^{(1)}_i \zeta^{(1)}_j e^{\zeta^{(1)} \cdot x} \; \zeta^{(2)}_k \zeta^{(2)}_l e^{\zeta^{(2)} \cdot x} \, dx \\
 =  \int_\Omega &\Big[ C_{1111}(\mathbf{i}s-t)^2(\mathbf{i}s+t)^2 + C_{1133}(\mathbf{i}s-t)^2(\mathbf{i}t-s)^2  \\
 + &C_{1313} (\mathbf{i}s-t)(\mathbf{i}t+s)(\mathbf{i}s+t)(\mathbf{i}t-s) + C_{1331}(\mathbf{i}s-t)(\mathbf{i}t+s)(\mathbf{i}t-s)(\mathbf{i}s+t) \\
+ &C_{3113}(\mathbf{i}t+s)(\mathbf{i}s-t)(\mathbf{i}s+t)(\mathbf{i}t-s)+ C_{3131} (\mathbf{i}t+s)(\mathbf{i}s-t)(\mathbf{i}t-s)(\mathbf{i}s+t) \\
 + &C_{3311} (\mathbf{i}t+s)^2 (\mathbf{i}s+t)^2 + C_{3333} (\mathbf{i}t+s)^2 (\mathbf{i}t-s)^2 \Big] e^{2\mathbf{i}(s,0,t)\cdot x} \, dx
\end{align*}
Combining the terms, one has
$$
0 = \int_\Omega (t-is)^2 (t+is)^2 [C_{1111} - 2C_{1133} + 4C_{1313} + C_{3333}] e^{2\mathbf{i}(s,0,t)\cdot x} \,dx.
$$
This means $\mathcal{F}\left[C_{1111} - 2C_{1133} + 4C_{1313} + C_{3333}\right](-2(s,0,t)) (t^2+s^2)^2 = 0$. A similar argument as in Step 1 shows $\mathcal{F}\left[\chi_\Omega (C_{1111} - 2C_{1133} + 4C_{1313} + C_{3333})\right](\xi) = 0$ for any $\xi \neq 0$, hence
\begin{equation} \label{eq:simpleCGO}
C_{1111} - 2C_{1133} + 4C_{1313} + C_{3333} = 0, \quad\quad\quad \text{ in } \Omega.
\end{equation}

\bigskip

\noindent\textbf{Step 3.}
We still take $u,v$ of the form
$$
u =  \zeta^{(1)} e^{\zeta^{(1)} \cdot x}, \quad\quad\quad 
v = \zeta^{(2)} e^{\zeta^{(2)} \cdot x},
$$
but with different phases $\zeta^{(1)},\,\zeta^{(2)}$. Set $d :=\sqrt{s^2+t^2}$ and $\beta:=\sqrt{\frac{r^2}{d^2}-1}$. The new phases to be used are
\begin{align*}
\zeta^{(1)} & = \mathbf{i}(s,0,t) + \mathbf{i}\beta (-t,0,s) + r(0,1,0) = (\mathbf{i}s-\mathbf{i}\beta t, r, \mathbf{i}t+\mathbf{i}\beta s), \\
\zeta^{(2)} & = \mathbf{i}(s,0,t) - \mathbf{i}\beta (-t,0,s) - r(0,1,0) = (\mathbf{i}s+\mathbf{i}\beta t,-r, \mathbf{i}t-\mathbf{i}\beta s). 
\end{align*}
It is easy to verify $\zeta^{(1)}\cdot \zeta^{(1)} = \zeta^{(2)}\cdot \zeta^{(2)} = 0$. This property again makes $\Delta u = \Delta v = 0$ and $\nabla \cdot u = \nabla \cdot v = 0$. Note that these new phases include the old ones: they coincide with \eqref{phase1} if one takes $r=0$ and $\beta=-\mathbf{i}$.

%
%
%
%
%

Using such $u,v$ in $\eqref{eq:intiden}$, we have
$$
0 =  \int_\Omega C_{ijkl} \zeta^{(1)}_i \zeta^{(1)}_j e^{\zeta^{(1)} \cdot x} \; \zeta^{(2)}_k \zeta^{(2)}_l e^{\zeta^{(2)} \cdot x} \, dx =: G_1 + G_2 + G_3,
$$
where
\[
\begin{split}
G_1 := & \int_\Omega \Big[ C_{1111} (\mathbf{i}s-\mathbf{i}\beta t)^2 (\mathbf{i}s+\mathbf{i}\beta t)^2 + C_{2222} (r)^2 (-r)^2 \\
&~~~\quad\quad\quad\quad\quad\quad~~~~+ C_{3333} (\mathbf{i}t+\mathbf{i}\beta s)^2 (\mathbf{i}t-\mathbf{i}\beta s)^2 \Big] e^{2\mathbf{i}(s,0,t)\cdot x}  \,dx \\
 = & \int_\Omega \left[ C_{1111} (s-\beta t)^2 (s+\beta t)^2+ C_{2222} r^4 + C_{3333} (t+\beta s)^2 (t-\beta s)^2 \right] e^{2\mathbf{i}(s,0,t)\cdot x}  \,dx, 
 \end{split}
 \]

\[
\begin{split}
G_2 := & \int_\Omega \left[ C_{1122} (\mathbf{i}s-\mathbf{i}\beta t)^2 (-r)^2 + C_{2211} (r)^2 (\mathbf{i}s+\mathbf{i}\beta t)^2 \right. \\
 &\quad\quad + C_{1133} (\mathbf{i}s-\mathbf{i}\beta t)^2 (\mathbf{i}t-\mathbf{i}\beta s)^2  + C_{3311} (\mathbf{i}t+\mathbf{i}\beta s)^2 (\mathbf{i}s+\mathbf{i}\beta t)^2 \\
 & \quad\quad\quad\quad\quad\quad\quad\quad\left. + C_{2233} (r)^2 (\mathbf{i}t-\mathbf{i}\beta s)^2 + C_{3322} (\mathbf{i}t+\mathbf{i}\beta s)^2 (-r)^2  \right] e^{2\mathbf{i}(s,0,t)\cdot x} \, dx \\
= & \int_\Omega \left[  -2 C_{1122} (s^2+\beta^2 t^2) r^2  \right. \\
 &\quad\quad + C_{1133} (s-\beta t)^2 (t-\beta s)^2  + C_{3311} (t+\beta s)^2 (s+\beta t)^2 \\
 & \quad\quad\left.  -2 C_{2233} (t^2+\beta^2 s^2) r^2  \right] e^{2\mathbf{i}(s,0,t)\cdot x} \, dx, 
 \end{split}
 \]
 \[
\begin{split}
G_3 := & \int_\Omega \Big[  4 C_{1212} (\mathbf{i}s-\mathbf{i}\beta t)(r)(\mathbf{i}s+\mathbf{i}\beta t)(-r)  \\
&\quad\quad\quad\quad\quad\quad+ 4 C_{1313} (\mathbf{i}s-\mathbf{i}\beta t) (\mathbf{i}t+\mathbf{i}\beta s) (\mathbf{i}s+\mathbf{i}\beta t) (\mathbf{i}t-\mathbf{i}\beta s)  \\
 &\quad\quad\quad\quad\quad\quad + 4 C_{2323} (r)(\mathbf{i}t+\mathbf{i}\beta s)(-r)(\mathbf{i}t-\mathbf{i}\beta s) \Big] e^{2\mathbf{i}(s,0,t)\cdot x} \, dx\quad\quad\quad\quad\quad\quad\quad \\
 = & \int_\Omega \Big[  4 C_{1212} (s^2 - \beta^2 t^2) r^2  +4 C_{1313} (s^2-\beta^2 t^2) (t^2 - \beta^2 s^2)\\
&\quad\quad\quad\quad\quad\quad\quad\quad\quad+  4 C_{2323} (t^2-\beta^2 s^2) r^2  \Big] e^{2\mathbf{i}(s,0,t)\cdot x} \, dx.
 \end{split}
 \]

We will analyze the asymptotic behavior as $r\rightarrow \infty$. Direct calculation (though tedious) shows
\begin{align*}
G_1 = & \int_\Omega \left(C_{1111} \frac{t^4}{d^4}+C_{2222}+C_{3333}\frac{s^4}{d^4}\right)r^4 e^{2\mathbf{i}(s,0,t)\cdot x}  \,dx+O(r^3), \\
G_2 = & \int_\Omega \left(-2C_{1122} \frac{t^2}{d^2}+2C_{1133}\frac{t^2s^2}{d^4}-2C_{2233}\frac{s^2}{d^2}\right)r^4 e^{2\mathbf{i}(s,0,t)\cdot x}  \,dx+O(r^3), \\
G_3 = & \int_\Omega \left(-4C_{1212} \frac{t^2}{d^2}+4C_{1313}\frac{t^2s^2}{d^4}-4C_{2323}\frac{s^2}{d^2}\right)r^4 e^{2\mathbf{i}(s,0,t)\cdot x}  \,dx+O(r^3),
\end{align*}
Equating the terms of order $O(r^4)$ yields
\begin{align*}
 \int_\Omega \Big[ & \frac{t^4}{d^4} C_{1111} + C_{2222} + \frac{s^4}{d^4} C_{3333} - \frac{2t^2}{d^2} C_{1122} + \frac{2 t^2 s^2}{d^4} C_{1133} \\
& \quad- \frac{2s^2}{d^2} C_{2233} - \frac{4 t^2}{d^2} C_{1212} + \frac{4 t^2 s^2}{d^4} C_{1313} - \frac{4 s^2}{d^2}  \Big] e^{2\mathbf{i}(s,0,t)\cdot x} \, dx=0.
\end{align*}
Using the linear relations in \eqref{eq:linrel} and $d^2=t^2+s^2$, one can eliminate $C_{2222}, C_{2233}, C_{1212}, C_{2323}$ and get
\begin{align*}
 \int_\Omega \Big[ \frac{s^4}{d^4} C_{1111}  - \frac{2 s^4}{d^4} C_{1133} - \frac{4 s^4}{d^4} C_{1313} + \frac{s^4}{d^4} C_{3333}  \Big] e^{2\mathbf{i}(s,0,t)\cdot x} \, dx=0.
\end{align*}
This, combined with the same argument used before, implies the following linear relation:
\begin{equation} \label{eq:r4}
 C_{1111}  - 2 C_{1133} - 4 C_{1313} +  C_{3333}=0 \quad\quad\quad \text{ in } \Omega. 
\end{equation}

%
%
%
%
%
%

Let us put the three pieces of information \eqref{eq:r5}\eqref{eq:simpleCGO}\eqref{eq:r4} together
\begin{align*}
C_{1111} & -   C_{1122} + 0\cdot C_{1133} +  2 C_{1313} + 0 \cdot C_{3333} = 0; \\
C_{1111} & + 0 \cdot C_{1122} - 2 C_{1133} + 4 C_{1313} + C_{3333} = 0;\\
 C_{1111} & + 0 \cdot C_{1122} - 2 C_{1133} - 4 C_{1313} +  C_{3333} = 0.
\end{align*}
We observe these combinations are linearly independent and thus can be used to eliminate $3$ independent components of $C$. In fact, solving this linear system yields 
\begin{align}
C_{1313} = & C_{1212} = \frac{1}{2}(C_{1111}-C_{1122})=0; \label{eq:knownrel1}\\
2C_{1133} = & C_{1122}+C_{3333} \label{eq:knownrel2}.
\end{align}
We are therefore left with only $2$ independent components, say $C_{1111}$ and $C_{1133}$. 

~\\
\noindent \textbf{The need for different solutions.}
The previous CGOs are not enough to determine the remaining independent components. To see this, we employ the known relations \eqref{eq:knownrel1}\eqref{eq:knownrel2} to simplify the integral identity $(\ref{eq:intiden})$, then
\begin{equation}\label{eq:e1}
\begin{split}
0 =  \int_\Omega &C_{1111} (\partial_1u_1\partial_1v_1+\partial_1u_1\partial_2v_2+\partial_2u_2\partial_1v_1+\partial_2u_2\partial_2v_2)\\
+& C_{1133}(\partial_1u_1\partial_3v_3+\partial_3u_3\partial_1v_1+\partial_2u_2\partial_3v_3+\partial_3u_3\partial_2v_2)\\
+& C_{3333}\partial_3u_3\partial_3v_3\, dx\\
= \int_\Omega &C_{1111} \nabla\cdot u\nabla\cdot v\\
 +& (C_{1133}-C_{1111})(\partial_1u_1\partial_3v_3+\partial_3u_3\partial_1v_1+\partial_2u_2\partial_3v_3+\partial_3u_3\partial_2v_2)\\
+& (C_{3333}-C_{1111})\partial_3u_3\partial_3v_3\, dx.\\
=\int_\Omega &C_{1111} \nabla\cdot u\nabla\cdot v\\
 +& (C_{1133}-C_{1111})(\partial_1u_1\partial_3v_3+\partial_3u_3\partial_1v_1+\partial_2u_2\partial_3v_3+\partial_3u_3\partial_2v_2)\\
+& 2(C_{1133}-C_{1111})\partial_3u_3\partial_3v_3\, dx.\\
=\int_\Omega &C_{1111} \nabla\cdot u\nabla\cdot v\\
 +& (C_{1133}-C_{1111})(\partial_1u_1\partial_3v_3+\partial_3u_3\partial_1v_1+\partial_2u_2\partial_3v_3+\partial_3u_3\partial_2v_2+2\partial_3u_3\partial_3v_3)\, dx.\\
 =\int_\Omega &C_{1111} \nabla\cdot u\nabla\cdot v+ (C_{1133}-C_{1111})(\nabla\cdot u\partial_3v_3+\partial_3u_3\nabla\cdot v)\, dx.
\end{split}
\end{equation}
All the solutions we have constructed have divergence zero, so they cannot give new information about the tensor $\mathbf{C}$. 

\begin{remark}
With only CGO solutions of divergence zero, one cannot even determine an isotropic perturbation from $\dot{\Lambda}_{\mathbf{C}^0}$ (cf. \cite{Ike1}). To see this, suppose $C_{ijkl}=\lambda\delta_{ij}\delta_{kl}+\mu(\delta_{ik}\delta_{jl}+\delta_{il}\delta_{jk})$, then \eqref{eq:intiden} reduces (with $\delta C_{ijkl}$ abbreviated as $C_{ijkl}$) to
\begin{equation}\label{Ale}
\int_{\Omega}2\mu\mathrm{Sym}(\nabla u):\mathrm{Sym}(\nabla v) + \lambda(\nabla \cdot u)(\nabla \cdot v)dx=0.
\end{equation}
Here $\mathrm{Sym}(\nabla u):=\frac{1}{2}(\nabla u+(\nabla u)^T)$ and $A:B = \sum^3_{i,j=1} A_{ij} B_{ij}$ for any $3\times 3$ matrices $A, B$. It is obvious that solutions with divergence zero cannot provide information about $\lambda$.
\end{remark}

%
\bigskip
\noindent \textbf{New type of solutions.}
This above analysis suggests the necessity to construct new solutions with non-vanishing divergence. We proceed to construct new CGO-type solutions with this property. 
They are of the form
\[
u=[(b\cdot x)\hat{\zeta}+c]e^{\zeta\cdot x}
\]
where $\zeta\in\mathbb{C}^3$ satisfies $\zeta \cdot \zeta = 0$, $\hat{\zeta}$ denotes $\frac{\zeta}{|\zeta|}$, and $b,c$ are constant vectors to be determined. This type of solutions can be constructed as in \cite{Ike1}.
The divergence of $u$ is
$$
\nabla\cdot u = \nabla\cdot\left([(b\cdot x)\hat{\zeta}+c]e^{\zeta\cdot x}\right) = [(b\cdot x)\hat{\zeta}\cdot\zeta+b\cdot\hat{\zeta}+c\cdot\zeta]e^{\zeta\cdot x} = (b\cdot\hat{\zeta}+c\cdot\zeta)e^{\zeta\cdot x},
$$
hence
\[
\nabla\nabla\cdot u=(b\cdot\hat{\zeta}+c\cdot\zeta)\zeta e^{\zeta\cdot x}.
\]
On the other hand, the gradient of $u$ is
\[
\nabla u=(b\otimes \hat{\zeta}+\zeta\otimes c+(b\cdot x)\zeta\otimes \hat{\zeta})e^{\zeta\cdot x}
\]
so $\Delta u$ can be computed:
$$
\Delta u = \nabla\cdot\nabla u
=[(\hat{\zeta}\cdot b)\zeta+(\zeta\cdot\zeta)c+(b\cdot x)(\zeta\cdot\zeta)\hat{\zeta}+(b\cdot \hat{\zeta})\zeta]e^{\zeta\cdot x}
=2(\hat{\zeta}\cdot b)\zeta e^{\zeta\cdot x}.
$$
We then have
\[
\begin{split}
\mu^0\Delta u + (\lambda^0+\mu^0)\nabla\nabla\cdot u=&[2\mu^0(\hat{\zeta}\cdot b)\zeta + (\lambda^0+\mu^0)(b\cdot\hat{\zeta}+c\cdot\zeta)\zeta]e^{\zeta\cdot x}\\
=&[(\lambda^0+\mu^0)c\cdot\zeta + (\lambda^0+3\mu^0)b\cdot\hat{\zeta}]\zeta e^{\zeta\cdot x}.
\end{split}
\]
Taking $b=(\lambda^0+\mu^0)\Re\hat{\zeta}$ and $c=-\frac{\lambda^0+3\mu^0}{|\zeta|}\Re\hat{\zeta}$ guarantees the right hand side is zero, making $u$ a solution to $(\ref{eq:elastic2})$.
%
Notice that with such $b,c$, the divergence of $u$ is
$$
\nabla\cdot u = (b\cdot\hat{\zeta} + c\cdot\zeta)e^{\zeta\cdot x}
= -2\mu^0\Re\hat{\zeta}\cdot\hat{\zeta}\,e^{\zeta\cdot x} = -\mu^0\,e^{\zeta\cdot x},
$$
which is non-vanishing since $\mu^0 >0$.

\bigskip

\noindent\textbf{Step 4.}
We take
$$
u =  \zeta^{(1)} e^{\zeta^{(1)} \cdot x},\quad\quad v=[(b\cdot x)\widehat{\zeta^{(2)}}+c]e^{\zeta^{(2)}\cdot x},
$$
with
\begin{equation*}
\begin{split}
\zeta^{(1)} & := \mathbf{i}(s,0,t) + \mathbf{i}\beta (-t,0,s) + r(0,1,0) = (\mathbf{i}s-\mathbf{i}\beta t, r, \mathbf{i}t+\mathbf{i}\beta s), \\
\zeta^{(2)} & := \mathbf{i}(s,0,t) - \mathbf{i}\beta (-t,0,s) - r(0,1,0) = (\mathbf{i}s+\mathbf{i}\beta t,-r, \mathbf{i}t-\mathbf{i}\beta s).
\end{split}
\end{equation*}
It has been verified that $\zeta^{(1)}\cdot \zeta^{(1)} = \zeta^{(2)}\cdot \zeta^{(2)} = 0$; moreover, $|\zeta^{(1)}|=|\zeta^{(2)}|=\sqrt{2}r$. Correspondingly, we take
$$
b= (\lambda^0+\mu^0)\Re \widehat{\zeta^{(2)}}=(0,-\frac{\lambda^0+\mu^0}{\sqrt{2}},0),\quad\quad\quad c=\frac{\lambda^0+3\mu^0}{\sqrt{2}r}(0,\frac{1}{\sqrt{2}},0).
$$
Substitute $u,v$ into $(\ref{eq:e1})$ and notice $\nabla\cdot v=-\mu^0\,e^{\zeta^{(2)}\cdot x}$, $\nabla\cdot u=0$. we have
\[
\begin{split}
0=-&\int_\Omega(C_{1133}-C_{1111}) \mu^0\zeta^{(1)}_3\zeta^{(1)}_3 e^{(\zeta^{(1)}+\zeta^{(2)})\cdot x}\, dx\\
=&\mu^0 \, \int_\Omega(C_{1133}-C_{1111}) (t+\beta s)^2 e^{2\mathbf{i} (s,0,t)\cdot x}\, dx \\
= &\mu^0\, \left[ \int_\Omega(C_{1133}-C_{1111}) \frac{s^2}{d^2} e^{2\mathbf{i} (s,0,t)\cdot x}\, dx\right]  r^2 + O(r)
\end{split}
\]
where the asymptotics is again when $r\rightarrow\infty$.
This implies 
\begin{equation}\label{identity4}
\frac{s^2}{t^2+s^2} \int_\Omega(C_{1133}-C_{1111}) e^{2\mathbf{i} (s,0,t)\cdot x}\, dx = 0
\end{equation}
In other words, $\mathcal{F}[\chi_\Omega (C_{1133}-C_{1111})](-2(s,0,t))=0$ when $s\neq 0$. Using the definition of transverse isotropy, one sees this is true in the entire $\mathbb{R}^3$ except for a plane (the one corresponding to $\{s=0\}$). Moreover, the Fourier transform is actually an analytic function since $\chi_\Omega (C_{1133}-C_{1111})$ is compactly supported. This forces $\mathcal{F}[\chi_\Omega (C_{1133}-C_{1111})]=0$ everywhere, so $C_{1133}=C_{1111}$ in $\Omega$.


\bigskip

\noindent\textbf{Step 5.}
Now we have $C_{1111}=C_{1133}=C_{3333}$, and $(\ref{eq:e1})$ becomes
\begin{equation}\label{eq:e2}
0 =  \int_\Omega C_{1111}\nabla\cdot u\nabla\cdot v\,dx=0.
\end{equation}
Take
$$
u=[(b^{(1)}\cdot x)\widehat{\zeta^{(1)}}+c^{(1)}]e^{\zeta^{(1)}\cdot x},\quad\quad\quad v=[(b^{(2)}\cdot x)\widehat{\zeta^{(2)}}+c^{(2)}]e^{\zeta^{(2)}\cdot x},
$$
with
\begin{equation*}
\begin{split}
\zeta^{(1)} & := \mathbf{i}(s,0,t) + (-t,0,s), \\ 
\zeta^{(2)} & := \mathbf{i}(s,0,t) - (-t,0,s),
\end{split}
\end{equation*}
$$
b^{(1)} = -b^{(2)} = \left(-\frac{t(\lambda^0+\mu^0)}{\sqrt{2}d},0,\frac{s(\lambda^0+\mu^0)}{\sqrt{2}d}\right), \quad\quad\quad 
$$
$$
c^{(1)} = -c^{(2)} = (\lambda^0+3\mu^0)\left(\frac{t}{2d^2},0,-\frac{s}{2d^2}\right).
$$
Substitue $u$ and $v$ into $(\ref{eq:e2})$
\[
\int_{\Omega}C_{1111}\left(\mu^0\right)^2 e^{2\mathbf{i}(s,0,t)\cdot x}\, dx=0.
\]
 Then we get $C_{1111}=0$. This completes the proof of the uniqueness of all parameters in $\mathbf{C}$.

\bigskip
\section{Non-zero frequency case: Proof of Theorem \ref{thm:maink}}

We present the proof of Theorem \ref{thm:maink} in this section. Notice that $\dot{\Lambda}_{\mathbf{C}^0,\rho^0,\omega}(\delta \mathbf{C},\delta\rho)=0$ implies
\begin{equation} \label{eq:intidenk}
\int_\Omega \delta C_{ijkl}(x) \, \partial_i u_j(x) \, \partial_k v_l(x) -\omega^2\delta\rho_{ik}u_iv_k\,dx=0,
\end{equation}
where $u,v$ satisfies \eqref{eq:elastic2k}.\\

\noindent \textbf{Step 1.}
For any $s^2+t^2>k_\mathbf{s}^{2}:=\omega^2c_\mathbf{s}^{-2}$, set 
\begin{equation}\label{phase1k}
\zeta^{(1)} := \mathbf{i}(s,0,t) + (-t,0,s)\sqrt{1-\frac{k^2_\mathbf{s}}{s^2+t^2}}, \quad\quad\quad 
\zeta^{(2)} := \mathbf{i}(s,0,t) - (-t,0,s)\sqrt{1-\frac{k^2_\mathbf{s}}{s^2+t^2}},
\end{equation}
where $c_\mathbf{s}^2:=\frac{\mu^0}{\rho^0}$ is the speed of \textit{S}-wave,
and
$$
a^{(1)} = a^{(2)} = a = (0,1,0).
$$ 
Take
$$
u = a^{(1)}e^{\zeta^{(1)} \cdot x},\quad\quad v = a^{(2)}e^{\zeta^{(2)} \cdot x}.
$$
It is easy to verify that  $\zeta^{(1)}\cdot \zeta^{(1)} = \zeta^{(2)}\cdot \zeta^{(2)} = -k^2_\mathbf{s}$, $a\cdot\zeta^{(j)}=0,\,j=1,2$ and $u,v$ solve the equations $(\ref{eq:elastic2k})$.
Denote $\mathfrak{K}:=\sqrt{1-\frac{k^2_\mathbf{s}}{s^2+t^2}}$, substitute $u,v$ into $(\ref{eq:intidenk})$, we have
\begin{equation}\label{A1k}
\begin{split}
0 = & \int_\Omega (C_{ijkl} a_i \zeta^{(1)}_j\; a_k \zeta^{(2)}_l -\omega^2 \rho_{ik}a_ia_k)e^{(\zeta^{(1)}+\zeta^{(2)}) \cdot x} \, dx \vspace{1ex}\\
= & \int_\Omega [C_{1212}(\mathbf{i}s-\mathfrak{K}t)(\mathbf{i}s+\mathfrak{K}t)+C_{1313}(\mathbf{i}t+\mathfrak{K}s)(\mathbf{i}t-\mathfrak{K}s)] e^{2\mathbf{i}(s,0,t)\cdot x}\\
&\quad\quad\quad\quad\quad\quad\quad \omega^2\rho_{11}e^{2\mathbf{i}(s,0,t)\cdot x}\, dx\\
= & \int_\Omega [(s^2+t^2)(-C_{1212}-C_{1313})+\frac{\omega^2c_\mathbf{s}^{-2}}{s^2+t^2}\left(t^2C_{1212}+s^2C_{1313}\right)-\omega^2\rho_{11}]e^{2\mathbf{i}(s,0,t)\cdot x}\, dx.
\end{split}
\end{equation}
If we have above identity at two different frequencies $\omega=\omega_1,\omega_2$, we can separate the two parts to obtain
\begin{equation}\label{eq:r11k}
\mathcal{F}[\chi_\Omega(C_{1212}+C_{1313})](2s,0,2t)=0
\end{equation}
and
\begin{equation} \label{eq:r2k}
\frac{c_\mathbf{s}^{-2}}{s^2+t^2}\left(t^2\mathcal{F}[ \chi_\Omega C_{1212}](2s,0,2t)+s^2\mathcal{F}[\chi_\Omega C_{1313}](2s,0,2t)\right)-\mathcal{F}[\chi_\Omega\rho_{11}](2s,0,2t)=0.
\end{equation}
By the transverse isotropy assumption, \eqref{eq:r11k} implies $\mathcal{F}[\chi_\Omega(C_{1212}+C_{1313})](\xi)=0$ for any $|\xi|\geq \sqrt{2}k_\mathbf{s}$. Then by the analyticity of the Fourier transform of compactly supported functions, we have
\begin{equation}\label{eq:r1k}
C_{1212}+C_{1313}=0\quad\quad\quad\text{in}~\Omega.
\end{equation}
This is exactly the same identity as what we got in Step 1 for the zero-frequency case.
\bigskip

\noindent\textbf{Step 2.}
The proof will be quite different from the zero-freqeuncy case from now on.
For $s^2+t^2>k_\mathbf{p}^{2}:=\omega^2c_\mathbf{p}^{-2}$, where $c_\mathbf{p}^2=\frac{\lambda^0+2\mu^0}{\rho^0}$, take $u,v$ of the form
$$
u =  \zeta^{(1)} e^{\zeta^{(1)} \cdot x}, \quad\quad\quad 
v = \zeta^{(2)} e^{\zeta^{(2)} \cdot x},
$$
with 
\begin{equation}\label{phase3k}
\begin{split}
\zeta^{(1)} & = \mathbf{i}(s,0,t) + \mathbf{i}\beta (-t,0,s) + r(0,1,0) = (\mathbf{i}s-\mathbf{i}\beta t, r, \mathbf{i}t+\mathbf{i}\beta s), \\
\zeta^{(2)} & = \mathbf{i}(s,0,t) - \mathbf{i}\beta (-t,0,s) - r(0,1,0) = (\mathbf{i}s+\mathbf{i}\beta t,-r, \mathbf{i}t-\mathbf{i}\beta s),
\end{split}
\end{equation}
where $d :=\sqrt{s^2+t^2}$ and $\beta:=\sqrt{\frac{r^2}{d^2}-1+\frac{k_\mathbf{p}^{2}}{d^2}}$. 
It is easy to verify $\zeta^{(1)}\cdot \zeta^{(1)} = \zeta^{(2)}\cdot \zeta^{(2)} =k_\mathbf{p}^2$.  Therefore the solutions $u,v$ constructed here is not divergence free. Substitute $u,v$ into $(\ref{eq:intidenk})$, we have
\begin{equation}\label{Aleeq1}
0 =  \int_\Omega C_{ijkl} \zeta^{(1)}_i \zeta^{(1)}_j e^{\zeta^{(1)} \cdot x} \; \zeta^{(2)}_k \zeta^{(2)}_l e^{\zeta^{(2)} \cdot x} \, dx =: \int_\Omega (\sum_{j=1}^9I_j+\omega^2 \sum_{k=1}^3J_k)e^{2\mathbf{i}(s,0,t)\cdot x} \, dx. 
\end{equation}
Here
\begin{equation*}
\begin{split}
I_1&=C_{1111} (s-\beta t)^2 (s+\beta t)^2,\\
I_2&=C_{2222} r^4,\\
I_3&=C_{3333} (t+\beta s)^2 (t-\beta s)^2,\\
I_4&=-2 C_{1122} (s^2+\beta^2 t^2) r^2,\\
I_5&=C_{1133} (s-\beta t)^2 (t-\beta s)^2  + C_{3311} (t+\beta s)^2 (s+\beta t)^2,\\
I_6&=-2 C_{2233} (t^2+\beta^2 s^2) r^2, \\
I_7&=4 C_{1212} (s^2 - \beta^2 t^2) r^2,\\
I_8&=4 C_{1313} (s^2-\beta^2 t^2) (t^2 - \beta^2 s^2),\\
I_9&=4 C_{2323} (t^2-\beta^2 s^2) r^2,
\end{split}
\end{equation*}
and
\begin{equation*}
\begin{split}
J_1&=\rho_{11}(s-\beta t)(s+\beta t),\\
J_2&=\rho_{11}r^2,\\
J_3&=\rho_{33}(t-\beta s)(t+\beta s).
\end{split}
\end{equation*}
We use the following asymptotics of $\beta$ and its powers in large $r$:
\begin{align*}
\beta & = \frac{r}{d}+\frac{(k_\mathbf{p}^2-d^2)r^{-1}}{2d}+O(r^{-3}), \\
 \beta^2&=\frac{r^2}{d^2}+\frac{k_\mathbf{p}^2-d^2}{d^2}+O(r^{-3}),\\
 \beta^4&=\frac{r^4}{d^4}+\frac{2r^2(k_\mathbf{p}^2-d^2)}{d^4}+\frac{(k_\mathbf{p}^2-d^2)^2}{d^4}+O(r^{-2}),
\end{align*}
and do some tedious calculation to obtain
\begin{align*}
I_1 & \sim r^4C_{1111}\frac{t^4}{d^4}+r^2C_{1111}\left(\frac{2(k_\mathbf{p}^2-d^2)t^4}{d^4}-\frac{2}{d^2}t^2s^2\right)+C_{1111}\left(\frac{d^2-k_\mathbf{p}^2}{d^2}t^2s^2+\frac{(d^2-k_\mathbf{p}^2)^2}{d^4}t^4+s^4\right) \\
I_2 & \sim r^4C_{1111} \\
I_3 & \sim r^4C_{3333}\frac{s^4}{d^4}+r^2C_{3333}\left(\frac{2(k_\mathbf{p}^2-d^2)s^4}{d^4}-\frac{2}{d^2}t^2s^2\right)+C_{3333}\left(\frac{d^2-k_\mathbf{p}^2}{d^2}t^2s^2+\frac{(d^2-k_\mathbf{p}^2)^2}{d^4}s^4+t^4\right) \\
I_4 & \sim r^4(-2C_{1122}\frac{t^2}{d^2})+r^22C_{1122}\left(\frac{d^2-k_\mathbf{p}^2}{d^2}t^2-s^2\right) \\
I_5 & \sim r^4(2C_{1133}\frac{t^2s^2}{d^4})+r^2\frac{2C_{1133}}{d^2}\left(s^4+t^4+2t^2s^2+2t^2s^2\frac{k_\mathbf{p}^2}{d^2}\right)\\
 & \quad\quad\quad\quad+2C_{1133}\left(t^2s^2\frac{(d^2-k_\mathbf{p}^2)^2}{d^4}-(s^4+t^4+4t^2s^2)\frac{d^2-k_\mathbf{p}^2}{d^2}+t^2s^2\right)
 \end{align*}
\begin{align*}
I_6 & \sim r^4\left(-2C_{2233}\frac{s^2}{d^2}\right)+r^22C_{2233}\left(\frac{d^2-k_\mathbf{p}^2}{d^2}s^2-t^2\right) \\
I_7 & \sim r^4\left(-4C_{1212}\frac{t^2}{d^2}\right)+r^24C_{1212}\left(\frac{d^2-k_\mathbf{p}^2}{d^2}t^2+s^2\right);\\
I_8 & \sim r^44C_{1313}\frac{t^2s^s}{d^4}+r^24C_{1313}\left(-2t^2s^2\frac{d^2-k_\mathbf{p}^2}{d^4}-\frac{t^4+s^4}{d^2}\right)\\
&\quad\quad\quad\quad+4C_{1313}\left(t^2s^2\frac{(d^2-k_\mathbf{p}^2)^2}{d^4}+\frac{d^2-k_\mathbf{p}^2}{d^2}(t^4+s^4)+s^2t^2\right);\\
I_9 & \sim r^4\left(-4C_{2323}\frac{s^2}{d^2}\right)+r^2 4C_{2323}\left(\frac{d^2-k_\mathbf{p}^2}{d^2}s^2+t^2\right).
\end{align*}
We also have
\begin{align*}
J_1 & \sim r^2(-\rho_{11})\frac{t^2}{d^2}+\rho_{11}\left(s^2+\frac{d^2-k_\mathbf{p}^2}{d^2}t^2\right); \\
J_2 & \sim r^2\rho_{11}; \\
J_3 & \sim r^2(-\rho_{33})\frac{s^2}{d^2}+\rho_{33}\left(t^2+\frac{d^2-k_\mathbf{p}^2}{d^2}s^2\right).
\end{align*}
The $O(r^4)$ terms in \eqref{Aleeq1} are exactly the same as in Step 3 for the zero-frequency case, which give
\begin{equation} \label{eq:r6k}
 C_{1111}  - 2 C_{1133} - 4 C_{1313} +  C_{3333}=0 \quad\quad\quad \text{ in } \Omega. 
\end{equation}
The coefficient of $r^2$ in the integrand is
\[
\begin{split}
\omega^2\Bigg[&c_\mathbf{p}^{-2}\Bigg(\frac{2C_{1111}t^4}{d^4}+\frac{2C_{3333}s^4}{d^4}-2C_{1122}\frac{t^2}{d^2}+4\frac{C_{1133}}{d^2}t^2s^2-2C_{2233}\frac{s^2}{d^2}\\
&-4C_{1212}\frac{t^2}{d^2}+8C_{1313}\frac{t^2s^2}{d^4}-4C_{2323}\frac{s^2}{d^2}\Bigg)-\rho_{11}\frac{t^2}{d^2}+\rho_{11}-\rho_{33}\frac{s^2}{d^2}\Bigg]\\
&+\Bigg[-C_{1111}\frac{2t^4}{d^2}-C_{1111}\frac{2t^2s^2}{d^2}-2C_{3333}\frac{s^4}{d^2}-C_{3333}\frac{2t^2s^2}{d^2}\\
&+2C_{1122}(t^2-s^2)+2C_{1133}\frac{(s^2+t^2)^2}{d^2}+2C_{2233}(s^2-t^2)+4C_{1212}(t^2+s^2)\\
&+4C_{1313}(-\frac{2t^2s^2}{d^2}-\frac{t^4+s^4}{d^2})+4C_{2323}(s^2+t^2)
\Bigg]
\end{split}
\]
Write this expression as $\omega^2A(x,s,t)+B(x,s,t)$, then two different frequency data separates the coefficients of $\omega^2$ and $\omega^0$ and gives
\begin{equation}\label{eq:r8k}
\int_{\Omega}A(x,s,t)e^{2\mathbf{i}(s,0,t)\cdot x} \, dx=\int_{\Omega}B(x,s,t)e^{2\mathbf{i}(s,0,t)\cdot x} \, dx=0.
\end{equation}
Using the relation \eqref{eq:linrel}, we get
$$
B(x,s,t)=\left(2 s^2 C_{1111} - 4 s^2 C_{1122} + 4 s^2 C_{1133} - 2 s^2 C_{3333}\right). 
$$
Then \eqref{eq:r8k} implies
\begin{equation}\label{eq:r9k}
C_{1111} - 2 C_{1122} + 2 C_{1133} - C_{3333}=0 \quad\quad\quad~~~\text{in}~\Omega.
\end{equation}

We summarize identities \eqref{eq:r1k}\eqref{eq:r6k}\eqref{eq:r9k} as
\begin{align*}
 1\cdot C_{1212} & + 0 \cdot( C_{1122}-2C_{1133}+C_{3333}) +  1\cdot C_{1313}=0, \\
 2\cdot C_{1212} & + 1 \cdot(C_{1122}-2C_{1133}+C_{3333}) -4\cdot C_{1313}=0,\\
  2\cdot C_{1212} & - 1 \cdot (C_{1122}-2C_{1133}+C_{3333}) + 0\cdot C_{1313}=0.
\end{align*}
By solving the above linear equations, we have $C_{1212}=C_{1313}=0$ and $2C_{1133} =  C_{1122}+C_{3333}$. Then \eqref{A1k} becomes $\mathcal{F}[\chi_\Omega\rho_{11}](2s,0,2t)=0$. Then we can get 
$$\rho_{11}=0\quad\quad\quad~~~\text{in}~\Omega.$$  

\bigskip
\noindent \textbf{Step 3.}
Take
$$
u = \vartheta^{(1)} e^{\zeta^{(1)} \cdot x}, \quad\quad\quad 
v =  \vartheta^{(2)} e^{\zeta^{(2)} \cdot x},
$$
with $\zeta^{(j)}$ defined in \eqref{phase1k} and 
\begin{equation}\label{amplitude1k}
\vartheta^{(1)} := \mathbf{i}(s,0,t)\sqrt{1-\frac{k_\mathbf{s}^{2}}{s^2+t^2}} + (-t,0,s), \quad\quad\quad 
\vartheta^{(2)} := \mathbf{i}(s,0,t)\sqrt{1-\frac{k_\mathbf{s}^{2}}{s^2+t^2}} - (-t,0,s).
\end{equation}
Notice that $\vartheta^{(1)}\cdot\zeta^{(1)}=\vartheta^{(2)}\cdot\zeta^{(2)}=0$. The solution $u,v$ used here is divergence free.

Inserting $u,v$ into the integral identity \eqref{eq:intidenk} and use the fact $C_{1313}=0$, $2C_{1133} =  C_{1122}+C_{3333}$ and $\rho_{11}=0$, we obtain
\begin{align*}
0 = \int_\Omega & C_{ijkl} \vartheta^{(1)}_i \zeta^{(1)}_j e^{\zeta^{(1)} \cdot x} \; \vartheta^{(2)}_k \zeta^{(2)}_l e^{\zeta^{(2)} \cdot x}-\omega^2\rho_{ik} \vartheta^{(1)}_ie^{\zeta^{(1)} \cdot x}\vartheta^{(2)}_ke^{\zeta^{(2)} \cdot x}, dx \vspace{1ex}\\
 =  \int_\Omega &\Big[ C_{1111}(\mathfrak{K}^2s^2+t^2)(\mathfrak{K}^2t^2+s^2) + C_{1133}(\mathbf{i}\mathfrak{K}s-t)(\mathbf{i}s-\mathfrak{K}t)(\mathbf{i}\mathfrak{K}t-s)(\mathbf{i}t-\mathfrak{K}s)   \\
 + &C_{3311} (\mathbf{i}\mathfrak{K}t+s)(\mathbf{i}t+\mathfrak{K}s)(\mathbf{i}\mathfrak{K}s+t)(\mathbf{i}s+\mathfrak{K}t) + C_{3333} (\mathfrak{K}^2s^2+t^2)(\mathfrak{K}^2t^2+s^2) \\
 +&\omega^2\rho_{33}(\mathfrak{K}^2t^2+s^2)\Big] e^{2\mathbf{i}(s,0,t)\cdot x} \, dx \vspace{1ex} \\
 = \int_\Omega &\Big[(C_{1111}-2C_{1133}+C_{3333})(s^2+\mathfrak{K}^2t^2)(t^2+\mathfrak{K}^2s^2)+\omega^2\rho_{33}(\mathfrak{K}^2t^2+s^2)\Big]e^{2\mathbf{i}(s,0,t)\cdot x} \, dx \vspace{1ex} \\
 =\int_\Omega &\omega^2[\rho_{33}(\mathfrak{K}^2t^2+s^2)]e^{2\mathbf{i}(s,0,t)\cdot x} \, dx  \vspace{1ex} \\
 =\int_\Omega &\omega^2[\rho_{33}(s^2+t^2)-\omega^2c_\mathbf{s}^{-2}\rho_{33}\frac{t^2}{s^2+t^2}]e^{2\mathbf{i}(s,0,t)\cdot x} \, dx
\end{align*}
Again using the above identity at two different frequencies, we end up with
$$\rho_{33}=0\quad\quad\quad\quad\text{in}~\Omega.$$

%

\bigskip
At this point we have recovered the density $\delta\rho$ and the same quantities for the elastic tensor $\delta\mathbf{C}$ as the zero-frequency case. Similar to \eqref{eq:e1} we have
\begin{equation}\label{eq:e1k}
\int_\Omega C_{1111} \nabla\cdot u\nabla\cdot v+ (C_{1133}-C_{1111})(\nabla\cdot u\partial_3v_3+\partial_3u_3\nabla\cdot v)\, dx=0.
\end{equation}

\bigskip
\noindent\textbf{Step 4.}
Take
$$
u =  \zeta^{(1)} e^{\zeta^{(1)} \cdot x}, \quad\quad\quad 
v = \vartheta^{(2)} e^{\zeta^{(2)} \cdot x},
$$
with $\zeta^{(1)},\zeta^{(2)}$ defined as \eqref{phase3k} and
\[
\vartheta^{(2)}=(\mathbf{i}s+\mathbf{i}\beta t,-\frac{r^2+k_\mathbf{p}^{2}}{r}, \mathbf{i}t-\mathbf{i}\beta s)
\]
Substitute into \eqref{eq:e1k}, and notice $\nabla\cdot u = -k_\mathbf{p}^{2}u$ and $\nabla\cdot v=0$. Then we obtain the following identity similar to \eqref{identity4} from the leading order terms in $r$:
$$
\frac{s^2}{t^2+s^2} \int_\Omega(C_{1133}-C_{1111}) e^{2\mathbf{i} (s,0,t)\cdot x}\, dx = 0,
$$
from which we get $C_{1133}=C_{1111}$.

\bigskip
\noindent\textbf{Step 5.}
For the last step, simply take
$$
u =  \zeta^{(1)} e^{\zeta^{(1)} \cdot x}, \quad\quad\quad 
v = \zeta^{(2)} e^{\zeta^{(2)} \cdot x},
$$
with $\zeta^{(1)},\zeta^{(2)}$ defined as \eqref{phase3k}. Substitute into 
\begin{equation}\label{eq:e2k}
0 =  \int_\Omega C_{1111}\nabla\cdot u\nabla\cdot v\,dx=0,
\end{equation}
and notice $\nabla\cdot u = -k_\mathbf{p}^{2}u$, $\nabla\cdot v = -k_\mathbf{p}^{2}v$. We can easily obtain $C_{1111}=0$, and conclude the proof of Theorem \ref{thm:maink}.

%
%
%
%
%
%
%
\bibliographystyle{abbrv}
\bibliography{biblio}

\begin{thebibliography}{10}

\bibitem{TI2}
T.~Alkhalifah.
\newblock Velocity analysis using nonhyperbolic moveout in transversely
  isotropic media.
\newblock {\em Geophysics}, 62(6):1839--1854, 1997.

\bibitem{TI3}
T.~Alkhalifah, I.~Tsvankin, K.~Larner, and J.~Toldi.
\newblock Velocity analysis and imaging in transversely isotropic media:
  {M}ethodology and a case study.
\newblock {\em The Leading Edge}, 15(5):371--378, 1996.

\bibitem{Ruiz2018uniqueness}
J.~Barceló, M.~Folch-Gabayet, S.~Pérez-Esteva, A.~Ruiz, and M.~Vilela.
\newblock Uniqueness for inverse elastic medium problems.
\newblock {\em SIAM Journal on Mathematical Analysis}, 50(4):3939--3962, 2018.

\bibitem{Bauchau2009book}
O.~A. Bauchau and J.~I. Craig.
\newblock {\em Structural Analysis with Applications to Aerospace Structures}.
\newblock Springer-Verlag, New York, 2009.

\bibitem{Beretta2017uniqueness}
E.~Beretta, M.~V. de~Hoop, E.~Francini, S.~Vessella, and J.~Zhai.
\newblock Uniqueness and lipschitz stability of an inverse boundary value
  problem for time-harmonic elastic waves.
\newblock {\em Inverse Problems}, 33(3):035013, 2017.

\bibitem{BFMRV}
E.~Beretta, E.~Francini, A.~Morassi, E.~Rosset, and S.~Vessella.
\newblock Lipschitz continuous dependence of piecewise constant {L}am{\'e}
  coefficients from boundary data: the case of non-flat interfaces.
\newblock {\em Inverse Problems}, 30(12):125005, 2014.

\bibitem{BFV}
E.~Beretta, E.~Francini, and S.~Vessella.
\newblock Uniqueness and {L}ipschitz stability for the identification of
  {L}am\'{e} parameters from boundary measurements.
\newblock {\em Inverse Probl. Imag.}, 8:611--644, 2014.

\bibitem{bhattacharyya2018local}
S.~Bhattacharyya.
\newblock Local uniqueness of the density from partial boundary data for
  isotropic elastodynamics.
\newblock {\em Inverse Problems}, 34(12):125001, 2018.

\bibitem{Calderon}
A.~P. Calder{\'o}n.
\newblock On an inverse boundary value problem.
\newblock {\em Computational \& Applied Mathematics}, 25(2-3):133--138, 2006.

\bibitem{CHN}
C.~I. C\^{a}rstea, N.~Honda, and G.~Nakamura.
\newblock Uniqueness in the inverse boundary value problem for piecewise
  homogeneous anisotropic elasticity.
\newblock {\em SIAM J. Math. Anal.}, 50(3):3291--3302, 2018.

\bibitem{ER}
G.~Eskin and J.~Ralston.
\newblock On the inverse boundary value problem for linear isotropic
  elasticity.
\newblock {\em Inverse Problems}, 18(3):907, 2002.

\bibitem{greenleaf}
A.~Greenleaf, M.~Lassas, and G.~Uhlmann.
\newblock On nonuniqueness for {C}alder{\'o}n's inverse problem.
\newblock {\em Mathematical Research Letters}, 10(5):685--693, 2003.

\bibitem{Ike1}
M.~Ikehata.
\newblock Inversion formulas for the linearized problem for an inverse boundary
  value problem in elastic prospection.
\newblock {\em SIAM Journal on Applied Mathematics}, 50(6):1635--1644, 1990.

\bibitem{Ike2}
M.~Ikehata.
\newblock The linearization of the {D}irichlet to {N}eumann map in anisotropic
  plate theory.
\newblock {\em Inverse Problems}, 11(1):165, 1995.

\bibitem{Ike4}
M.~Ikehata.
\newblock The linearization of the {D}irichlet-to-{N}eumann map in the
  anisotropic {K}irchhoff-{L}ove plate theory.
\newblock {\em SIAM Journal on Applied Mathematics}, 56(5):1329--1352, 1996.

\bibitem{Ike3}
M.~Ikehata.
\newblock A relationship between two {D}irichlet to {N}eumann maps in
  anisotrpoic elastic plate theory.
\newblock {\em Journal of Inverse and Ill-Posed Problems}, 4(3):233--244, 1996.

\bibitem{IUY}
O.~Y. Imanuvilov, G.~Uhlmann, and M.~Yamamoto.
\newblock On uniqueness of {L}am\'{e} coefficients from partial {C}auchy data
  in three dimensions.
\newblock {\em Inverse Problems}, 28(12):125002, 2012.

\bibitem{IY}
O.~Y. Imanuvilov and M.~Yamamoto.
\newblock Global uniqueness in inverse boundary value problems for the
  {N}avier--{S}tokes equations and {L}am{\'e} system in two dimensions.
\newblock {\em Inverse Problems}, 31(3):035004, 2015.

\bibitem{TI1}
D.~Kumar, M.~K. Sen, and R.~J. Ferguson.
\newblock Traveltime calculation and prestack depth migration in tilted
  transversely isotropic media.
\newblock {\em Geophysics}, 69(1):37--44, 2004.

\bibitem{LN}
Y.-H. Lin and G.~Nakamura.
\newblock Boundary determination of the {L}am\'{e} moduli for the isotropic
  elasticity system.
\newblock {\em Inverse Problems}, 33:125004, 2017.

\bibitem{TI4}
J.-P. Montagner and H.-C. Nataf.
\newblock A simple method for inverting the azimuthal anisotropy of surface
  waves.
\newblock {\em Journal of Geophysical Research: Solid Earth}, 91(B1):511--520,
  1986.

\bibitem{nachman1988reconstructions}
A.~I. Nachman.
\newblock Reconstructions from boundary measurements.
\newblock {\em Annals of Mathematics}, 128(3):531--576, 1988.

\bibitem{NU3}
G.~Nakamura and G.~Uhlmann.
\newblock Inverse problems at the boundary for an elastic medium.
\newblock {\em SIAM J. Math. Anal.}, 26:263--279, 1995.

\bibitem{NTU}
G.~Nakamura and G.~Uhlmann.
\newblock Layer stripping for a transversely isotropic elastic medium.
\newblock {\em SIAM J. Appl. Math.}, 59:1879--1891, 1999.

\bibitem{NU2}
G.~Nakamura and G.~Uhlmann.
\newblock Erratum: Global uniqueness for an inverse boundary value problem
  arising in elasticity.
\newblock {\em Invent. Math.}, 152:205--207, 2003.
\newblock Erratum to Invent. Math., 118, 457-474, (1994).

\bibitem{rachele2000inverse}
L.~Rachele.
\newblock An inverse problem in elastodynamics: uniqueness of the wave speeds
  in the interior.
\newblock {\em Journal of Differential Equations}, 162(2):300--325, 2000.

\bibitem{rachele2003uniqueness}
L.~Rachele.
\newblock Uniqueness of the density in an inverse problem for isotropic
  elastodynamics.
\newblock {\em Transactions of the American Mathematical Society},
  355(12):4781--4806, 2003.

\bibitem{Reddy2002book}
J.~N. Reddy.
\newblock {\em Energy principles and variational methods in applied mechanics}.
\newblock John Wiley \& Sons, New York, second edition, 2002.

\bibitem{stefanov2017local}
P.~Stefanov, G.~Uhlmann, and A.~Vasy.
\newblock Local recovery of the compressional and shear speeds from the
  hyperbolic {DN} map.
\newblock {\em Inverse Problems}, 34(1):014003, 2017.

\bibitem{SU1}
J.~Sylvester and G.~Uhlmann.
\newblock A global uniqueness theorem for an inverse boundary value problem.
\newblock {\em Ann. of Math.}, 125:153--169, 1987.

\end{thebibliography}

\end{document}